\documentclass{amsart}
\newtheorem{theorem}{Theorem}[section]
\newtheorem{lemma}[theorem]{Lemma}
\newtheorem{corollary}[theorem]{Corollary}

\theoremstyle{definition}

\theoremstyle{remark}
\newtheorem{remark}[theorem]{Remark}

\numberwithin{equation}{section}

\def \real{{\mathbb R}}

\def \integer{{\mathbb Z}}

\def \Id {{\text{Id}\, }}
\def \Lip {{\text{Lip}\, }}

\def\FF{{\mathcal F}}

\def\HH{{\mathcal H}}

\def\LL{{\mathcal L}}

\def\today {\ifcase\month\or January \or February \or March \or
April \or May \or June
\or July \or August \or September \or October \or November \or December
\fi
\number\day~\number\year}

\begin{document}
\title[Semi-flows without  finite Markov partitions]
{Exponential decay of correlations for surface 
semi-flows without finite Markov partitions}

\author{Viviane Baladi}
\address
{CNRS, UMR 7586, 
Institut de Math\'e\-ma\-ti\-ques de Jussieu, 75251 Paris, FRANCE}

\email{baladi@math.jussieu.fr}

\author{Brigitte Vall\'ee}

\address{CNRS, GREYC, Universit\'e de Caen, 14032 Caen, FRANCE}

\email{brigitte.vallee@info.unicaen.fr}

\date{November 2003}

\begin{abstract}
We extend Dolgopyat's bounds on iterated transfer
operators to suspensions of 
interval maps with infinitely many intervals of monotonicity.
\end{abstract}


\thanks{Work supported
in part by two CNRS MATHSTIC grants. We thank A. Hachemi for
a careful reading of a previous version of the paper.}

\maketitle

%%%%%%%%%%%%%%%%%%%%%%%%%%%%%%%%%%

\section{Statement of results}

Let $0<c_1 < \ldots < c_m<c_{m+1} <\ldots <1$ be a finite or countable
partition of $I=[0,1]$ into subintervals,
and let $T:I\to I$ be so that $T|_{(c_m, c_{m+1})}$
is  $C^2$ and extends to a 
homeomorphism from $[c_m, c_{m+1}]$ to $I$. 
We assume that $T$ is piecewise uniformly expanding:
there are
$C \ge 1$ and $\hat \rho < 1$  so that $|h(x)-h(y)|\le C \hat\rho^n
|x-y|$ for every inverse branch $h$ of $T^n$ and all $n$.
Let $\HH$ be the  set of inverse branches
$h:I\to [c_m, c_{m+1}]$ of $T$. 
We  suppose
(Renyi's condition) that there is  $K>0$ so that every 
$h\in \HH$ satisfies $|h''|\le K |h'|$.
Let $r:I \to \real_+$ be so that $r|_{(c_m, c_{m+1})}$ is
$C^1$, and $\inf r > 0$. Assume  that
there is $\sigma_0<0$ so that
$\sum_{h \in \HH}\sup \exp(-\sigma (r\circ h)) |h'| <\infty$ for all
$\sigma >  \sigma_0$, and that 
$ |r' \circ h| \cdot |h'| \le K$ 
for all $h\in \HH$.  
For $n\ge 1$, write
$r^{(n)}(x)=\sum_{k=0}^{n-1}r(T^k)(x)$.

We study the transfer operators, indexed by $s=\sigma+it$,
$$
L_s f(x) = \sum_{T(y)=x }e^{-s r(y)}  \frac{ f(y)}{ |T'(y)|}=
\sum_{h \in \HH} e^{-s r(h(x))}  
 |h'(x)| \cdot (f\circ h)(x) \, ,
$$
acting on  $C^1(I)$,
with norm $\|f\|=\sup|f|+\sup|f'|$. 
Note that $L_s=\LL_{s+h_{top}}$, where
$\LL_s$ is the transfer operator associated to the suspension  semi-flow
on the branched surface $\{(x,s)\in I\times \real_+\mid s\le r(x)\}/\sim$,
with $(x,r(x))\sim (T(x), 0)$,  
defined by $\phi^t(x,s)=(x,s+t)$, and $h_{top}$ is the topological
entropy of $\phi^t$. See e.g. \cite{P}.

Finally, the following assumption
is a translation of  Dolgopyat's ``uniform nonintegrability
of foliations'' condition (see  \cite{A, P, PS} for  formulations
closer to ours): we say that the pair
$(T,r)$  satisfies {\it UNI} if there exist 
$D>0$ and  $n_0 \ge 1$ such that, for every
integer $n \ge n_{0}\ge 1$, there are
two elements $h$, $k$ of the set $\HH_n$ of inverse
branches of $T^{n}$ so that
the function on $I$ defined by
$\psi_{h,k}(x):=r^{(n)} (h (x))-r^{(n)}(k(x))$
satisfies
$\inf |\psi_{h,k}'| \ge D$.
(See also Remark~\ref{remUNI}.)

To state our main result, we use the  equivalent norms
$\|f\|_{1,t}=\sup|f|+ \frac{\sup|f'|}{ |t|}$,  for
$|t|\ge \epsilon_0 > 0$, on $C^1(I)$:

\begin{theorem}\label{T1}
Let $T$ and $r$ satisfy the assumptions above (in particular {\it UNI} for $D$). 
Then there is $A \ge n_0$ and $\gamma < 1$
so that
for all $\sigma$ close enough to $0$, 
all  $|t|\ge \max(2\pi/D, 4 )$, 
and all $n \ge A \log |t|$,  we have
$
\| L_s^{ n} \|_{1,t}\le \gamma^n$.
\end{theorem}

Theorem \ref{T1} was proved by Dolgopyat \cite{Do} when $\HH$ is
finite. In \cite{BV}, we considered the special case when $T(x)=\{1/x\}$
(or analogues of the Gauss map) and $r=\log |T'|$, working with a 
different version of {\it UNI}, adapted to ``algebraic'' situations.
Note that the present {\it UNI} assumption also holds in the
setting of \cite{BV}: if $h\in \HH_n$ is a linear fraction
$(ax+b)/(cx+d)$ then $h''(x)/h'(x)=-2c/(cx+d)$ so that
$|\psi'_{h,\hat h}(x)|=|2[c/(cx+d)-\hat c/(\hat cx+\hat d)]|$. 
Write $\FF_n$ for the $n$th Fibonacci
number and $\widehat \FF_n$ for the sequence $0$, $1$, 
$\widehat \FF_n= 2\widehat \FF_{n-1}+\widehat \FF_{n-2}$. For $h$
and $\hat h$ in $\HH_n$ associated to
the sequence of digits $1,1,\ldots, 1$, and $2,2,\ldots, 2$, we get 
$a=\FF_{n-2}$, $b=c=\FF_{n-1}$, and $d=\FF_{n}$, while
$\hat a=\widehat\FF_{n-2}$, $\hat b=\hat c=\widehat \FF_{n-1}$, and
 $\hat d=\widehat \FF_{n}$. We conclude by using
$\lim_{n\to \infty} \FF_n/\FF_{n-1}=(1+\sqrt 5)/2$ and
$\lim_{n\to \infty} \widehat\FF_n/\widehat\FF_{n-1}=(1+\sqrt 8)/2$.

\smallskip

From Theorem \ref{T1}, one easily gets (see e.g. \cite{BV}):

\begin{corollary}\label{C2}
For every  $0 < \alpha < 1$ there is $t_0$ so that
for all  $|t|>t_0$ and 
$\sigma$ close to $0$, we have
$
\|(\Id-L_s)^{-1}\|_{1,t} \le  |t|^\alpha  
$.
\end{corollary}

Theorem \ref{T1}   implies \cite[section 4]{P}  exponential decay of correlations for
$C^ 1$ observables and the absolutely continuous invariant probability
(SRB) measure of the semi-flow $\phi^t$. 
We hope this will be a useful step towards proving
exponential decay of correlations for (continuous-time)
planar Sinai billiards, using \cite{Y}. (For the moment, only open 
continuous-time billiards
have been considered \cite{St}, they admit finite Markov sections.)
See Remark \ref{Gib} for extensions to
other Gibbs states.

\section{Proof of Theorem \ref{T1} }

We basically follow Dolgopyat's proof,
as detailed in \cite{P}, \cite{PS}, and  \cite{A}.
A key point is the {\it Federer property} of any absolutely
continuous measure $\nu$ with continuous
density bounded from above and from
below: There are $C, C'>0$ so that
if $I$, $J$ are adjacent intervals with $|I|\ge |J|/C$
then $\nu(I)\ge \nu(J)/C'$. To exploit this information
when considering $L_\sigma$ for $\sigma\ne 0$,
the arguments in \cite{Do} (e.g. last lines of p. 367) and 
\cite{A} (e.g. first lines of p. 43) use  that there
is $\alpha_\sigma \to 1$ when $\sigma \to 0$ so that 
$\widetilde L_\sigma f(x) \le \alpha_\sigma \widetilde L_0 f(x)$,
for  the normalised operators in (\ref{tA.6}) and positive $f$.
The above inequality uses that there are finitely many branches
and is for example not true for the Gauss map.
To bypass this problem, we  exploit carefully the Cauchy-Schwartz decomposition
in Lemma \ref{A.8} below  (see also 
%\cite{A}, \S 5.2 and 
\cite{BV}, Lemma 2).

\remark\label{Gib}
Beware that {\it even when there are finitely many
branches,} the Federer property  
is {\it not} true  for arbitrary Gibbs measures $\nu$, in particular
the measures $\nu_\sigma$ introduced below for $\sigma\ne 0$, contrary
to what is stated in \cite[Proposition 7]{Do}; 
\cite[Lemma 6]{P}; and \cite[Lemma 4]{PS}.  (Proposition 7 of
\cite{P} is true  e.g. if $T$ is a $C^2$ circle map,
and if $r$ is $C^1$ on the circle, and not only piecewise $C^1$.
For a counterexample, take $T(x)=2x$ modulo $1$ with 
$\exp(r)\equiv 3$
on $[0,1/2]$ and $\exp(r)\equiv 3/2$ on $(1/2, 1]$,  and
consider the intervals of size $1/2^n$ to the right
and to the left of $1/2$. By adding
$\epsilon \sin (2\pi x)$ to $r$, this example can probably be made
to satisfy the {\it UNI} condition \cite[p. 537]{P}.)
When there are finitely many branches,
the Federer property  {\it does} hold  \cite{Do'}
for Gibbs measures and ``most'' adjacent intervals from the partitions in  
\cite{Do}, 
\cite{P},  \cite{PS}: This is enough e.g. to recover the results
in \cite{Do}, in particular Theorem~1. 
When $\HH$ is infinite, the situation is more complicated
but we expect that Theorem~\ref{T1} will also hold for more
general transfer operators 
$L_{s,g} f(x) = \sum_{T(y)=x }e^{-s r(y)} g(y)  f(y)$ associated
to suitable positive $g$.
\endremark

{\bf Preliminary steps}

\smallskip

Fix from now on 
$\hat\rho  < \rho<1$.
The inverse branches of $T^n$
satisfy $|h''|\le \bar K |h'|$ for all $n$ and the distorsion constant $\bar K=K/(1-\rho)$.
Similarly for $(r^{(n)})'\circ h$.
As a consequence it is easy to prove that, for every $n\ge 1$, and each pair
$h$, $k$ in $\HH_n$, the function
$\psi_{h,k}(x)=r^{(n)}\circ h-r^{(n)}\circ k $ satisfies
$\sup |\psi_{h,k}'|=
\le 2 \bar K $.
We next recall  spectral properties of the $L_s$
(see e.g.  \cite{BV} and references therein).
Let $\sigma> \sigma_0$ be real. 
The essential spectral radius $\lambda^e_\sigma$ of
$L_\sigma$ is strictly smaller than its spectral radius
$\lambda_\sigma$ (in fact $\lambda^e_\sigma \le \rho \lambda_\sigma$).
Since $T$ is topologically
mixing, the operator $L_\sigma$ has a unique (simple) eigenvalue
$\lambda_\sigma$
of maximal modulus, for a strictly positive $C^1$ eigenfunction
$f_\sigma$, the rest of the spectrum is in the disc
of radius $\tau_\sigma \lambda_\sigma$ for some $\tau_\sigma < 1$. 
The eigenvector $\mu_\sigma$ of
$L_\sigma^*$ for $\lambda_\sigma$ is Lebesgue measure
for $\sigma=0$, and  for all $\sigma> \sigma_0$
a positive Radon measure $\mu_\sigma$.
We may  assume $\mu_\sigma(1)=1$ and
$\mu_\sigma(f_\sigma)=1$ so that $\nu_\sigma=f_\sigma \mu_\sigma$
is a probability measure.  
Note that $L_\sigma:C^1(I)\to C^1(I)$ depends continuously
on $\sigma$,
so that $\lambda_\sigma^{\pm1}$, $\tau_\sigma$, $f_\sigma^{\pm1}$, and $f'_\sigma$
depend continuously on $\sigma$ (and therefore satisfy
uniform bounds in any compact subset $\Sigma\subset (\sigma_0, \infty)$).
Also, $\sigma \mapsto \lambda_\sigma$ is a nonincreasing function.
Finally, the spectral radius of $L_{\sigma+it}$ is not larger than $\lambda_\sigma$
and its essential spectral radius is not larger than $\rho \lambda_\sigma$
for all $t\in \real$.

It will be convenient to work with the normalised operators
\begin{equation}
\label{tA.6}
\widetilde L_s( f) = \frac {L_s  (f_\sigma \cdot f)}
{\lambda_\sigma f_\sigma}\, , \quad s=\sigma+it \, .
\end{equation}
If $s=\sigma > \sigma_0$, the operator $\widetilde L_\sigma$ acting on $C^1(I)$
has spectral radius $1$, essential spectral
radius $ \le \rho$,  and fixes the constant function $\equiv 1$.
Clearly $\widetilde L_\sigma^*$ preserves the
probability measure $\nu_\sigma=f_\sigma \cdot \mu_\sigma$.
Our starting point is a Lasota-Yorke inequality:

\begin{lemma} \label{A.3} {\bf (Lasota-Yorke)}
For every compact $\Sigma\subset (\sigma_0, \infty)$, there is 
a constant 
$C =C(\Sigma,\bar K)>0$, so that for  all $n\ge 1$,
all $s\in \Sigma$ and all $f\in C^1(I)$:
\begin{equation}
\label{tA.7}
| (\widetilde  L^n_s f)'  (x) | \le   
C(\Sigma, \bar K) |s|  \cdot\widetilde  L^n_\sigma (|f |)(x) +
 \rho^n   \cdot \widetilde  L^n_\sigma(|f'| )(x)\, . 
\end{equation}
\end{lemma}

Lemma \ref{A.3} indicates that we should concentrate on the $\sup$
component of the norm, which will be estimated by the $\LL^2(d\nu_0)$
norm 
(see the beginning of the proof of Theorem~\ref{T1}
below). Indeed, the crucial estimate (Lemma ~\ref{A.8}) will show
exponential decay of iterates of the operator in
the $\LL^2(d\nu_0)$ norm.
These $\LL^2(d\nu_0)$ integrals are oscillatory integrals in disguise: because 
of the weights $\exp(-s (r\circ h))$ in $L_s$, the integrand is the absolute value
of a sum of  complex numbers with strong phase variations for
large $|t|$ ({\it UNI} is crucial here). We shall exhibit enough cancellations
in the terms, via the key Lemma~\ref{A.5}.

\begin{proof} 
The Leibniz sum for the derivative of each term 
$
\exp(-s(r^{(n)}\circ h))
 |h'|  \cdot\frac {1}{\lambda_\sigma f_\sigma } \cdot (f_\sigma f) \circ h 
$
forming $(\widetilde  L^n_s (f))'$ ($h\in \HH_n$) contains four terms.
We can bound the first for all $s$  using our ``distortion''
assumption on $r$ since
$|s| |(r^{(n)})'\circ h| |h'|  e^{-s)(r\circ h)} \le
 | s| \bar K  e^{-s)(r\circ h)}$.
The second one is controlled by the Renyi assumption on $T$.
Compactness of $\Sigma$ and continuity of
$\sigma\mapsto \lambda_\sigma$ and $\sigma\mapsto f_\sigma$ imply
$\sup_{\sigma \in \Sigma}|f'_\sigma |< \infty$  and
$\inf_{\sigma \in \Sigma}f_\sigma > 0$,
so that the third term may be controlled by
$
 \frac{|f'_\sigma| }{ \lambda_\sigma f_\sigma^2}\le C_\Sigma
\frac {1 }{ \lambda_\sigma f_\sigma}
$
for some $C_\Sigma > 0$.
Finally the last term can be estimated using
$$
 | (f_\sigma\cdot f) ' \circ h |  |h'| \le \rho^n  [ | f_\sigma' \cdot  f|  \circ h 
+   (f_\sigma \cdot |f'|)  \circ h ] \, .
$$
We can ensure
$
\bar K |s| + 2C_\Sigma + 2 \rho C_\Sigma \le C(\Sigma, \bar K) |s|
$  
(for fixed $\Sigma$, if $|s|$ is large, i.e., if $|t|$ is large
enough, then  $C(\Sigma, \bar K)$ is close to $\bar K$).
\end{proof}

We next state and prove an elementary 
lemma about complex numbers with almost opposite phases.
Note that  $2/3 <(\sqrt 7-1)/2< 1$.

\begin{lemma} \label{A.4} {\bf  (Calculus lemma)}
For each $\eta \in[ (\sqrt 7-1)/2, 1)$ and every pair of complex numbers,
$z_1=r_1 \exp(i\theta_1)$ and  $z_2=r_2\exp(i\theta_2)$,   
\begin{equation}
\label{tA.8}
 \cos(\theta_1-\theta_2) \le 1/2 \Rightarrow
|z_1+z_2|\le \max( \eta r_1 + r_2, r_1 + \eta r_2)\, .
\end{equation}
\end{lemma}

\begin{proof} 
Up to exchanging $z_1$ and $z_2$,
we can  suppose that $r_1 \le r_2$ so that $\eta r_1+r_2 \ge r_1+\eta r_2$.
Our assumption on $\cos(\theta_1-\theta_2)$ implies  
$$
|z_1+z_2|^2=r_1^2+r_2^2+ 2 r_1r_2 \cos(\theta_1-\theta_2)\le
r_1^2+r_2^2+r_1r_2\, .
$$
Since $(\eta r_1 + r_2)^2 = \eta^2 r_1^2 +r_2^2 +2 \eta r_1 r_2$,
we must show
$
r_1^2 (1-\eta^2) + 2r_1 r_2 (1/2-\eta)\le 0 
$.

Since   $\eta-1/2\ge 1-\eta^2 \ge 0$
(use $3\ge \sqrt 7\ge 2$), we get
$$
r_1^2 (1-\eta^2) + 2r_1 r_2 (1/2-\eta)\le r_1^2(\eta-1/2)+ r_1r_2(1/2-\eta)
\le r_1 (\eta-1/2) (r_1-r_2) \le 0 \, .
$$
\end{proof}

{\bf Preparatory lemmas in view of $\LL^2$ contraction}
\smallskip

In the next lemma, we combine  {\it UNI}  and Lemma \ref{A.4} to obtain 
cancellation-type estimates on terms appearing when applying iterates of $\widetilde L_\sigma$
to a  suitable pair $(u,v)$ of initial test functions in $C^1(I)$. We first introduce 
the  ``cone'' condition
that $(u,v)$ must satisfy: there are $C > 0$ and $t\in \real$ so that
\begin{equation}
\label{tA.9}
u > 0 \, , \quad 0\le |v| \le u\, , \quad
\max( |u'(x)|, |v'(x)|) \le 2 C |t| u(x)\, . 
\end{equation}

\begin{lemma}\label {A.5} {\bf (Exhibiting cancellations)}
Assume that {\it UNI} holds for $D$ and
$n_{0}$.
Then, for all $C >0$,  there are $n_1\ge n_{0}$,
$\delta > 0$  and $\Delta>0$, so that for any 
$|t| > 2\pi/D$, and all  $u, v \in C^1(I)$ 
satisfying (\ref{tA.9}) for $C$ and $t$, we have the following:

Fix $n \ge n_1$, and let $h, k \in \HH_n$ be the branches from {\it UNI.}
For every $x_0\in I$, there is $x_1\in I$ with $|x_0-x_1|< \Delta/|t|$, so that the function
$$
F(x):=
e^{-(\sigma+it) r^{(n)}(h(x))} |h'(x)|((v \cdot f_\sigma )\circ h)  (x)+
e^{-(\sigma+it)r^{(n)}(k(x))} |k'(x)]( (v \cdot f_\sigma )\circ k)(x)
$$
satisfies for  all  $x$ s.t. $|x -x_1| < \delta/|t|$,  all  $\sigma> \sigma_0$,
and all $\eta > (\sqrt 7 -1)/2$
\begin{equation}
\label{tA.10}
\begin{split}
&|F(x)| \le \max \bigl [\cr
&\eta  e^{-\sigma r^{(n)}(h(x))}|h'(x)| (u \cdot f_\sigma )\circ h)  (x)+
e^{-\sigma r^{(n)}(k(x))}|k'(x)|( (u \cdot f_\sigma )\circ k)(x),  \cr
&e^{-\sigma r^{(n)}(h(x))}|h'(x)|(u \cdot f_\sigma )\circ h)  (x)+
\eta e^{-\sigma r^{(n)}(k(x))}|k'(x)|( (u \cdot f_\sigma )\circ k)(x) )  \bigr ] \, .\cr
\end{split}
\end{equation}
When the maximum in (\ref{tA.10}) is attained by the expresssion
where the $\eta$ factor is attached to 
branch $h$ we say we are ``in case $h$,'' and otherwise ``in
case $k$.'' 
\end{lemma}

It follows from the proof  that $n_1\ge n_{0}$ so that $3 \times16 C \rho^{n_1}< 1/24$ 
works.
In the application of Lemma ~\ref{A.5} in Lemma~\ref{A.7} we require
$C\ge C(\Sigma, \bar K)$ from Lemma~ \ref{A.3}.

\begin{proof} 
Let us fix $x_0 \in I$. Assume first (this case does
not require {\it UNI}) that
\begin{equation}
\label{tA.11}
|v ( h(x_0))| \le \max(u(h(x_0))/2, u(k(x_0))/2 )\, . 
\end{equation}
Let us suppose the maximum is realised for   $u\circ h$ (the other case is
symmetric). Then it is easy to see that for any $\epsilon > 0$, if  
$|x-x_0| < \delta_1/|t|$,
with
$
\delta_1 (2C \rho^{n_{0}}) = \epsilon 
$,
we have
$
\exp(-\epsilon)  \le  \frac {u(h(x)) }{ u(h(x_0))}  \le \exp (\epsilon)
$.
(Use $\exp  ( \log u(h(x))-\log u(h(x_0)) ) \,  dx 
\le\exp \int_{h(x)}^{h(x_0)} 
|(\log u(y))'| \, dy$,
the assumed bound on $|u'|/u$ from (\ref{tA.9}) and
$n\ge n_{0}$.) 

To prove (\ref{tA.10}), it is then enough to check that
$|x-x_0| < \delta_1/|t|$ implies $|v (h(x) )| \le  \eta' u(h(x_0))$ for some
$\eta'>2/3$ with $\eta' \exp(\epsilon) \le \eta$: indeed, we would then
have
$|v (h(x) )| \le  \eta' \exp(\epsilon) u(h(x))\le \eta u(h(x))$  whenever  
$|x-x_0| < \delta_1/|t|$,
so that (\ref{tA.10}) would hold. 
Assume for a contradiction that no such $\eta'$
exists, i.e. for each $2/3 < \eta' \le \eta \exp(-\epsilon)$
there is $x_1$ with $|x_1-x_0|\le \delta_1/|t|$ and $|v( h(x_1) )| \ge  \eta' (u(h(x_0)))$, 
so that (use (\ref{tA.11}))
$
|v ( h(x_0)) -v(h(x_1))| \ge (\eta' -1/2) u(h(x_0))
$.
On the other hand, (\ref{tA.9}) and the choice
of $\epsilon$  imply that there is
$y$ with $|y-x_0|\le \delta_1/|t|$ so that
$$
|v ( h(x_0)) -v(h(x_1)) |\le u(h(y))  2C |t| \rho^{n_{0}}\frac  {\delta_1}{ |t|}
\le u(h(x_0)) e^ \epsilon 2C  \rho^{n_{0}}  \delta_1
= u(h(x_0)) \epsilon e^\epsilon\, ,
$$
a contradiction if $\epsilon \exp(\epsilon) < 1/6$.
This ends the easy case, where we can take $x_1=x_0$
(i.e. $\Delta_1=0$) and $\delta_1=\epsilon/(2C \rho^{n_{0}})$
for  small  (independently of $u$, $v$, $C$, etc.) $\epsilon>0$.
(The dependence of $\delta_1$ on $C$ can be removed by taking large enough
$n$.)

Let us now move to the more interesting situation when
\begin{equation}
\label{tA.12}
|v ( h(x_0))| > \max(u(h(x_0))/2, u(k(x_0))/2 )\, . 
\end{equation}
We shall use {\it UNI} to show that we are in a position to apply Lemma \ref{A.4}
to the sum forming $F(x)$, for $x$ in an $\delta_2/|t|$-interval around
a point $x_1$ which is $\Delta_2/|t|$ close to $x_0$.
Since $f_\sigma$ is real and positive, the difference  $\theta(x)$
between the argument of  the two terms of $F(x)$ can be decomposed as
$
\theta(x)= t\psi_{h,k}(x) + \arg (v(h(x))) -\arg (v(k(x)))
$.

Let us first show the claim by assuming that we found
$\delta_2$, $\Delta_2$ so that $\cos \theta(x)\le 1/2$,
for all $x$ with $|x-x_1|\le \delta_2 /|t|$, for some $x_1$ 
with $|x_1-x_0|<\Delta_2/|t|$,
leaving the (nontrivial) proof of this fact for the end.
We have $r_1(x)=e^{-\sigma r^{(n)}(h(x))}|h'(x)|( |v| \cdot f_\sigma )( h  (x))$ and
$r_2(x)=e^{-\sigma r^{(n)}(k(x))} |k'(x)|( |v |\cdot f_\sigma )( k(x))$.
Fix $x$ with $|x-x_1|\le \delta_2 /|t|$, and assume
(the other case is analogous) that $r_1(x)\le r_2(x)$.
Lemma \ref{A.4} then yields the claim:
\begin{equation*}
\begin{split}
|F(x)| &\le
\eta e^{-\sigma r^{(n)}(h(x))}  |h'(x)|(|v| \cdot f_\sigma )(h  (x))
+e^{-\sigma r^{(n)}(k(x))} |k'(x)|( |v| \cdot f_\sigma )( k(x))\cr
&\le
\eta e^{-\sigma r^{(n)}(h(x))}  |h'(x)| (u \cdot f_\sigma )(h  (x))
+e^{-\sigma r^{(n)}(k(x))} |k'(x)|( u \cdot f_\sigma )( k(x))\, .
\end{split}
\end{equation*}

It remains to prove that $\cos\theta(x)\le 1/2$ for $x$ as above and
some $\delta_2$, $\Delta_2$. For this,
the following consequence of (\ref{tA.12})
and (\ref{tA.9}) will be helpful: for all
$y, z$ with $|z-x_0| \le  |y-x_0| \le \xi/|t$
\begin{equation}
\begin{split}
\label{tA.13}
|v(h(y))| &\ge |v(h(x_0))| - |v(h(x_0))-v(h(y))|\cr
&\ge
u(h(x_0))/2 - \rho^{n_{0}}\frac {\xi}{ |t|} 2 C |t| u(h(z))\cr
&\ge  
u(h(x_0)) (1/2 - \exp(\epsilon) \rho^{n_{0}} \xi 2 C  )
\ge u(h(x_0))/4\, .
\end{split}
\end{equation}
Next observe that, because of (\ref{tA.9}, \ref{tA.12}),
$
V(x)=\arg (v(h(x)) -\arg (v(k(x))
$
does not vary too much around $x_0$.  More precisely:
\begin{equation}
\label{tA.14}
\begin{split}
|V(x)-V(x_0)| &=
| \log[ v(h(x)) /v(k(x))] - \log  [v(h(x_0)) / v(k(x_0)) ]\cr
&\le |\log [\frac {v(h(x)) }{ v(h(x_0))} ]| + |
 \log [\frac {v(k(x_0))}{ v(k(x))}  ] |\, ,
\end{split}
\end{equation}
and, if $|x-x_0|\le \xi/|t|$,
\begin{equation*}
\begin{split}
|\log [ \frac {v(h(x)) }{ v(h(x_0))} ]| 
&\le
 |h(x)-h(x_0)|   \frac{|v'(h(y))|}{ |v(h(y))|}\le
\rho^{n}  \frac{\xi}{ |t|} 8 C |t| e^\epsilon \frac {u(h(x_0)) }{  u(h(x_0))}
\le  \xi 8C e^ \epsilon\rho^{n}\, .
\end{split}
\end{equation*}
(We used $|y-x_0| \le |x-x_0|$  and (\ref{tA.13}, \ref{tA.12}).)
We may control $ | \log [\frac{v(k(x_0))}{ v(k(x))} ]|$, mutatis mutandis, and we have
for $|x-x_0|< \xi/|t|$:
\begin{equation}
\label{tA.15}
|V(x)-V(x_0)|\le  \xi 16C \exp(\epsilon)\rho^{n}\, .
\end{equation}

Recall that we have to show $\cos \theta(x)\le 1/2$
in a suitable interval. We first find $x_1$ with $|x_1-x_0| < \Delta_2/|t|$
such that  $|\theta(x_1)-\pi|\le \pi/24$.
For this, we use  {\it UNI} which ensures that, since
$t (\psi(z)-\psi(x_0))=t (z-x_0) \psi'(y)$ for $y\in I$,  if
$
\Delta_2=2\pi/D
$,
then
$
\{ t (\psi(z)-\psi(x_0))  \mod 2\pi \mid  |z-x_0|\le \Delta_2/|t| \}  =[0, 2\pi )
$.
(We use here $|t| > 2\pi/D$.)
In particular there is $z=x_1$ so that $t (\psi(x_1)-\psi(x_0)) =\pi-\theta(x_0)$
($\hbox{mod}\,  2\pi$).
Applying (\ref{tA.15}) to $x=x_1$, $\xi=\Delta_2$, we find
\begin{equation}
\label{tA.16}
\begin{split}
|\theta(x_1)-\pi|&=|\theta(x_0) + t(\psi(x_1)-\psi(x_0))+(V(x_1)-V(x_0))-\pi|\cr
&\le |V(x_1)-V(x_0)|< \Delta_2 16 C \exp(\epsilon)\rho^{n}< \pi/24 \, ,
\end{split}
\end{equation}
if $n$ is large enough  (depending on $C$ and, via $\Delta_2$, on  $D$).

To conclude, we apply (\ref{tA.15}) and the ``distorsion''
upper bound, using $|x-x_0| < |x-x_1|+|x_1-x_0|
< (\delta_2+\Delta_2)/|t|$ and $|x-x_1|<\Delta_2/|t|$ to get, if
$n$ is large enough (depending on $C$ and $D$)
and
$0<\delta_2\le \Delta_2$ is small enough (depending on $\bar K$):
\begin{equation}
\begin{split}
|\theta(x)-\pi|&\le \pi/24 + |\theta(x)-\theta(x_1)|\cr
&\le \pi/24 +| t| |\psi(x)-\psi(x_1)| + |V(x)-V(x_0)|+|V(x_1)-V(x_0)|\cr
&\le \pi/24+  2\bar K  |t| \frac{\delta_2 }{ |t|} 
+  16C \exp(\epsilon)\rho^{n} D |t| \frac{\delta_2+2\Delta_2 }{ |t| }
< \pi/12\, . 
\end{split}
\end{equation}

Taking $\delta=\min(\delta_1,\delta_2)$ and $\Delta=\Delta_2$,
we have proved the lemma.
\end{proof}

\remark\label{remUNI}
If we 
replace {\it UNI} by the assumption that there exist  $D>0$,   $n_{0}$,
and  two inverse branches $h$ and $k$ of $T^{n_{0}}$ so that
$\inf |\psi_{h,k}'|\ge D$,
then
for every $n\ge n_{0}$ there are 
$\hat h$, $\hat k \in \HH_n$ so that 
$\inf |\psi_{\hat h, \hat k}'|
\ge \rho^{n-n_{0}} D $. 
(Take $\hat h=h\circ \ell$, $\hat k=k\circ \ell$,
for $\ell \in \HH_{n-n_{0}}$ and observe that $\psi_{\hat h, \hat k} (x)=\psi_{h,k}(\ell(x))$.)
However, this is not enough.
In  (\ref{tA.16}) we would get
(in view of the definition of $\Delta_2$)
$
\frac
{2\pi}{ D\rho^{n-n_{0}} } 8C \exp(\epsilon)\rho^n=
\frac {16C}{ D}  \exp(\epsilon)\rho^{n_{0}}
$,
which is independent of $n$ and not necessarily smaller than $\pi/24$.
(The constant $16$ can be reduced, but not below $1$.)
Unfortunately, the strategy presented on p. 545 of \cite{P} seems to suffer from
the same problem.
\endremark

\noindent The following consequence of Lemma~\ref{A.5}
 will be instrumental towards
Lemma~\ref{A.8}:

\begin{corollary} \label{A.6}
Let $T$ satisfy {\it UNI} for $D$. Let $C >0$
and let $n_1=n_1(C)$, $\delta=\delta(C)$, $\Delta=\Delta(C)$ be given by Lemma ~\ref{A.5}.
Fix $n\ge n_1$, let $h, k\in \HH_n$ come from {\it UNI,} and
let $\rho_{n,C}=\min (\min |h'|,\min |k' |)$ (we have $0< \rho_{n,C}\le \rho^n$).

Then for every $|t|> 2\pi/D$,  every $u$, $v\in C^1(I)$
satisfying (\ref{tA.9}) for $C$ and $|t|$, and
each $\eta > (\sqrt 7 -1)/2$ (recall Lemma \ref{A.4}), there are:

\item $\bullet$
a finite set of (disjoint) intervals $[a_j, b_{j+1}] =I_j\subset I$,
$j=0, \dots, N-1$, with $|I_j| \ge \delta/|t|$,
$a_0 \le   \Delta/|t|$, and  $b_N \ge 1- \Delta/|t|$;
also, setting $J_j=[b_j, a_j]$, we have
$0 <|J_j| \le 2\Delta/|t|$;  to each $I_j$  is associated $type(I_j)
\in \{h, k\}$; we write $\hat I_j$ for the middle third interval of $I_j$;
\item $\bullet$
a function
$\chi=\chi(u,v,n,\eta) \in C^1(I)$ (in particular $\chi$ depends on $C$, $|t|$)
so that: 
$$
\begin{cases}
&\eta \le \chi(x) \le 1\, , \quad  |\chi'| \le
\frac {3(1-\eta) }{ \rho_{n,C} \delta} |t| \, ,\cr
&  type(I_j)=h \hbox{ and } x\in \hat I_j \Rightarrow  \chi_{h(x)}=\eta \, , \cr
&type(I_j)=k \hbox{ and } x \in \hat I_j  \Rightarrow \chi_{k(x)}= \eta \, , \cr
&\chi(y)<1 \Rightarrow   [y=h(x) \, , x\in I_j \, , type(I_j)=h]
\, \hbox{or} \,  [ y=k(x) \, , x\in I_j \, , type(I_j)=k ]   .
\end{cases}
$$

\noindent Finally,  for $s=\sigma+it$ with $\sigma> \sigma_0$, we have
$
|\widetilde L^n_s (v) (x) | \le \widetilde L_\sigma^n (\chi u) (x)\, , \, \, 
\forall x \in I 
$.
\end{corollary}

Note that  $\sup  |\chi'|/|t|$ can be made arbitrarily small by taking
$\eta<1$ close to $1$, once $C$ and $h$, $k$, $n$ are fixed.
To exploit Corollary~\ref{A.6}, we shall use the following:

\begin{lemma} \label{A.7} {\bf  (Invariance of  ``cone condition'')}
Let $T$ satisfy {\it UNI} for $D$ and fix $\Sigma$ a compact subset
of $(\sigma_0, \infty)$.  Let $C(\Sigma,\bar K)$ be from Lemma~\ref{A.3}  and fix $C >1$ so that:
$
C \ge C(\Sigma,\bar K)
\cdot \max(1, \max_{\sigma \in \Sigma} |\sigma| D/(2\pi)) 
$.

Then, there is $n_2\ge n_1$ ($n_1$ from Lemma~\ref{A.5})
so that for every large enough $|t|> 2\pi/D$, each
$u$, $v$, satisfying (\ref{tA.9}) for $C$ and $t$, 
and all $n\ge n_2$, taking  $\eta =\eta(n)< 1$ close enough to $1$,
and $\chi=\chi(u,v,\eta)$  from Corollary~\ref{A.6}, the pair
$
\hat u=\widetilde L^n_\sigma (\chi u)$,
$\hat v=\widetilde L^n_s (v)
$,
satisfies (\ref{tA.9}), for the same $|t|$ and $C$,
and for
all $s=\sigma+it$ with $\sigma\in \Sigma$.
\end{lemma}

\begin{proof} 
Corollary \ref{A.6} says that  
$
|\hat v (x)|=|\widetilde L^n_s (v) (x) | \le \widetilde L_\sigma^n (\chi u) (x)=
\hat u(x)$ for all $x\in I$.
We also have $\inf \hat u > 0$ since $\inf(\chi u) > 0$ and $\widetilde L_\sigma$
preserves the cone of strictly positive functions.
To check the condition on $\max (|\hat u'|, |\hat v'|)$ we shall (finally!) invoke the
Lasota-Yorke inequality from Lemma~ \ref{A.3}
(recalling also that $\widetilde L_\sigma$ is normalised so that
$\sup \widetilde L_\sigma |f|\le \sup |f|$). We first consider $\hat u'$
and get, using $|u'|\le 2C |t| u$,
$\chi\ge \eta$ and $|\chi'|\le 1$ ($\eta=\eta(C,n)$ is close to $1$):
\begin{equation*}
\begin{split}
\left | \frac{d}{ dx}\widetilde L_\sigma^n (\chi u) (x)\right |
&\le C(\Sigma, \bar K) \sigma \widetilde L_\sigma^n (\chi u) (x)
+ \rho^n \widetilde L_\sigma^n (|\chi' u+\chi u'|) (x)\cr
&\le 
\biggl ( C(\Sigma, \bar K) |t| + \rho^n (\frac{1}{ \eta}+2C|t|)\biggr )
\widetilde L_\sigma^n (\chi u) (x) 
\le  2 C |t| \widetilde L_\sigma^n (\chi u) (x)\, ,
\end{split}
\end{equation*}
if $n \ge n_2\ge n_1$ and $C \ge C(\Sigma,\bar K)$.

The computation for $|\hat v'|$ is similar:
\begin{equation}
\begin{split}
\left | \frac{d}{ dx} \widetilde L^n_s (v) (x) \right |
&\le  C(\Sigma, \bar K) |s| \widetilde L_\sigma^n (|v|) (x)
+ \rho^n \widetilde L_\sigma^n (|v'|) (x)\cr
&\le \frac{C(\Sigma, \bar K) |s| + 2C |t| \rho^n }{ \eta} \widetilde L_\sigma^n (\chi u) (x) 
\le  2 C |t| \widetilde L_\sigma^n (\chi u) (x)\, ,
\end{split}
\end{equation}
if $n \ge n_2\ge n_1$ and $C(\Sigma, \bar K) |s|\le C|t|$.
\end{proof}

{\bf Proof of the $\LL^2$ contraction and proof of Theorem~\ref{T1}}
\smallskip

We shall see below that 
the case $\sup |f'| > 2 C |t| \sup |f| $ 
 is easy. We next prove
the key ``$\LL^2$ contraction  lemma''
(see \cite[Lemma 4]{Do})  to handle the other  case:

\begin{lemma} \label{A.8} {\bf   ($\LL^2(\nu_1)$ contraction)}
Assume {\it UNI.} Let $\Sigma$, $C$,     $n\ge n_2$,
$|t|>2\pi/D$, be as in Lemma~\ref{A.7}.
There is $\beta < 1$ so that 
for all $\sigma$  close enough to $0$,
and for all $0\ne f \in C^1$ with
$\sup |f'| \le 2 C |t| \sup |f|$,
\begin{equation}
\label{tA.18}
\int |\widetilde L_{\sigma+it}^{m n} f |^2 \, d\nu_0
< \beta^m \sup |f|^2\, , \forall m \ge 1 , .
\end{equation}
\end{lemma}

\begin{proof} 
Recall $\eta<1$ was taken close to $1$ in Lemma \ref{A.7}.
For $s=\sigma+it$ with $\sigma \in \Sigma$, define a sequence
of pairs $(u_m, v_m)$, $m\ge 0$, of functions in $C^1(I)$:
\begin{equation*}
\begin{split}
&u_0\equiv 1\, , \, v_0=\frac {f}{ \sup |f|}\, , \,  \chi_0=\chi_{u_0, v_0, n}\, , 
\quad u_m = \widetilde L^{n}_\sigma(\chi_{m-1} u_{m-1}) \, ,\quad 
v_m = \widetilde L^{n}_s(v_{m-1}) \, .
\end{split}
\end{equation*}
Lemma~\ref{A.7} implies that all $(u_m, v_m)$ satisfy
(\ref{tA.9}) for $C$, $t$, and all $m$. 
(Note also that $u_m\le 1$ for all $m$.)
In particular,
$|\widetilde L_s^{m n}  (f/\sup |f| )|=|v_m|\le u_m$,
and to prove the lemma, it is enough to show that there is $\beta_1 < 1$,
so that 
$
\int u^2_{m+1}\, d \nu_0 \le  \beta_1 \int u^2_{m}\, d \nu_0\, 
$ for all $m$
(note that $\int u^2_0 \, d \nu_0=1$).

The definition of $u_{m+1}$ and the Cauchy-Schwartz inequality imply
\begin{equation*}
\begin{split}
\lambda_\sigma^{2n} f_\sigma^2(x) u_{m+1}^2(x)
&=
\biggl ( \sum_{\ell \in \HH_{n}} e^{-\sigma r^{(n)}(\ell(x))} |\ell'(x)|
(\chi_m \cdot f_\sigma \cdot u_m ) (\ell(x))
\biggr )^2\cr
&\le  \max_I \frac {f_\sigma }{ f_0}
\sum_{\ell \in \HH_{n}}  |\ell'(x)| (f_0 \cdot u_m^2 ) (\ell(x))\cr
&\qquad\qquad\quad\cdot
 \max_I \frac {f_\sigma }{ f_{2\sigma}}
\sum_{\ell \in \HH_{n}} e^{-2\sigma r^{(n)}(\ell(x))}
 |\ell'(x)| (\chi_m^2 \cdot f_{2\sigma} ) (\ell(x))
\, .
\end{split}
\end{equation*}
Now, if $x\in \hat I_j$ for $\chi_m$, of type $h$, say (type $k$
is similar), we get
\begin{equation*}
\begin{split}
\frac
{1}{ \lambda^n_{2\sigma} f_{2\sigma}(x)}&
\sum_{\ell \in \HH_{n}} 
e^{-2\sigma r^{(n)}\ell(x)}
 |\ell'(x)| (\chi_m^2 \cdot f_{2\sigma} ) (\ell(x))\cr
&\le
1 - (1-\eta^2)e^{-2\sigma r^{(n)}(h(x))} |h'(x)|
\frac { f_{2\sigma}  (h(x))}{ \lambda^n_{2\sigma} f_{2\sigma}(x)}
\le 1-\epsilon (1-\eta^2) =\eta' < 1
\end{split}
\end{equation*}
(we used $e^{-2\sigma r^{(n)}(h(x))}
|h'(x)|  f_{2\sigma}  (h(x))/( \lambda^n_{2\sigma} f_{2\sigma}(x))
\ge \epsilon > 0$
if $n$ and $h$ are fixed; obviously, $\eta'$ depends on $n$).
Denote
$$
\xi(\sigma,n)=
\frac{\lambda_{2\sigma}^n f_0(x) f_{2\sigma}(x)}{ \lambda_\sigma^{2n}  f_\sigma^2(x)}
\cdot \max_I \frac{f_\sigma }{ f_{2\sigma}}\cdot  \max_I \frac{f_\sigma }{ f_0}\, .
$$
We showed that for some $\eta' <1$ and all  $x \in \cup_j \hat I_j$ (recall $\lambda_0=1$)
$$
u_{m+1}^2(x)\le \eta'\xi(\sigma, n)  \widetilde L_0^{n} (u_m^2) (x) \, .
$$
If $x\notin \cup_j \hat I_j$, the Cauchy-Schwartz inequality  just gives, since $\chi_m\le1$,
$$
u_{m+1}^2(x)\le \xi(\sigma,n) \widetilde  L_0^{n} (u_m^2) (x) \, .
$$
We claim that there is $\hat \delta$, independent of
$m$, $n$, and $t$, so that if $\hat J_j$ is the union of the rightmost third of $I_j$, $J_j$, and
the leftmost third of $\hat I_{j+1}$, then
\begin{equation}
\label{tA.20}
\int_{\hat I_j} \widetilde L_0^{n}(u^2_m)\,  d\nu_0 \ge 
\hat \delta  \cdot \int_{\hat J_j} \widetilde L_0^{n}(u^2_m )\,  d\nu_0 \, .
\end{equation}
We finish the proof assuming (\ref{tA.20}): 
if $\hat \delta (\beta_2-\eta') \ge (1-\beta_2)$
(e.g. $\beta_2=:\frac {1+\eta'\hat \delta}{ 1+\hat \delta} <1$), 
\begin{equation}
\label{tA.21}
\begin{split}
\int_I u_{m+1}^2 \, d\nu_0
&\le \xi(\sigma, n)\sum_j \biggl ( \eta' \int_{\hat I_j}  \widetilde L_0^{n}(u_{m}^2) 
\, d\nu_0+
\int_{\hat J_j} \widetilde L_0^{n}(u_{m}^2 )\, d\nu_0 \biggr ) \cr
&\le  \xi(\sigma, n)
\beta_2 \biggl ( \sum_j \int_{\hat I_j}  \widetilde L_0^{n}(u_{m}^2) 
\, d\nu_0+
\int_{\hat J_j} \widetilde L_0^{n}(u_{m}^2 )\, d\nu_0 \biggr ) \cr
&= \xi(\sigma, n)
\beta_2 \int_I \widetilde  L_0^{n}(u_{m}^2 )\, d\nu_0 
= \xi(\sigma, n)
\beta_2 \int_I  u_{m}^2 \, d\nu_0 \, .
\end{split}
\end{equation}
(In the last line we used that the dual of   $\widetilde L_0^{n}$
leaves $\nu_0$ fixed.)
By taking $\sigma$ sufficiently close to $0$ (depending on $n$, which is fixed)
we can assume that $ \xi(\sigma, n)\cdot
\beta_2 <1$. 

It remains to show (\ref{tA.20}). It suffices to prove that 
$
\int_{\hat I_j} w^2 \,  d\nu_0 \ge \hat \delta \int_{\hat J_j} w^2 \,  
d\nu_0 $
for all $C^1$  functions $w$ with 
$\sup w\le 1$ and $|w'(z)|\le 2 C |t| w(z)$ (recall Lemma ~\ref{A.7}
and use Lemma~\ref{A.3} and $\widetilde \LL_0 1 \equiv 1$).
Note that such  $w$ satisfy, for all $x \in \hat I_j$,
$y\in \hat J_j$:
\begin{equation*}
\frac {w^2(y)}{ w^2(x)} =\exp 2(\log w(x)-\log w(y))=
\exp 2 \int_x^y (w'/w) (z) \, dz \le \exp( 4C (2\Delta+\delta))\, .
\end{equation*}
Applying the above inequality, and making use of the Federer
property (for intervals with length-ratio $3\Delta$), of $\nu_0$ which has density
$f_0$ (bounded from above and from below) with respect to Lebesgue measure,  we find
\begin{equation*}
\begin{split}
\int_{\hat I_j} w^2 \,  d\nu_0 &
\ge \nu_0(\hat I_j) \min_{\hat I_j} w^2
\ge  \tilde \delta  e^{-( 4C (2\Delta+\delta))} \nu_0(\hat J_j)  
\max_{\hat J_j} w^2
\ge \hat \delta \int_{\hat J_j} w^2 \,  d\nu_0 \, .
\end{split}
\end{equation*}
\end{proof}

We are finally ready to prove the theorem:

\begin{proof} 
Since there is $B$ so that ($\lambda_\sigma$ is semisimple)
$
\|L_s^{n} \|_{1,t}\le B \lambda_\sigma^{n} \|\widetilde L_s^{n} \|_{1,t}$
for all $n \ge 1$,
and since $\lambda_0=1$ and
$\sigma$ is in a neighbourhood of $0$,
 it is enough to show that there is  $\widetilde A$ 
and $\tilde  \gamma < 1$ so that
$
 \|\widetilde L_s^{n} \|_{1,t}\le \tilde \gamma^n
$,
for $n\ge \widetilde A \log |t|$. Clearly, this will follow from
the existence of $n_4$ and $\widehat A$ so that
$\|\widetilde L_s^{n_4 m} \|_{1,t}\le  \tilde \gamma^{n_4 m}$
for all 
$m\ge \widehat A \log |t|$
(write  $n=q n_4 + r$, with $q, r\in \integer^+$ and
$0 \le r < n_4$, and use
$\|\widetilde L_s^{q} \|_{1,t}\le \widetilde B$).

Let  (see Lemma \ref{A.7})
$
C=\max(3/2,C(\Sigma,\bar K)
\cdot \max(1,  D/(2\pi))) 
$, 
and  let $n_2$ be given by Lemma \ref{A.7}. Let $n_3\ge n_2$ be so that  $\rho^{n_3} < 1/4$.

Let us first deal with the easy case $\sup |f'| \ge 2 C |t| \sup |f| $.
Setting $\gamma_1=\max ((2C|t|)^{-1},  \rho^{n_3} +3/4) < 1$, we have
$
\sup |\widetilde L_s^{n_3 } f|\le \sup |f|\le \frac{1}{ 2C |t|} \sup |f'| \le 
\gamma_1\|f\|_{1,t} 
$,
and, by Lemma \ref{A.3},
\begin{equation}
\begin{split}
\frac { |(\widetilde L_s^{n_3} f)'| }{ |t|} 
&\le C(\Sigma, \bar K)\frac { |s|}{ |t|} \sup |f|+
\frac { \rho^{n_3 }}{ |t|}  \sup |f'| \cr
&\le \biggl (\frac { \sqrt {(\max |\sigma|^2 + |t|^2)}}{2 |t|}+  \rho^{n_3 }\biggr )  
\frac {\sup |f'| }{ |t|}
\le  \gamma_1 \|f\|_{1,t}  \, . \cr
\end{split}
\end{equation}

If  $\sup |g'| < 2 C |t| \sup |g| $, then
the function $g^2$ satisfies (\ref{tA.9})
for  $2C \max(1, \sup |g|)$ for which
Lemmas \ref{A.7}, \ref{A.8} hold. Note also that a slight modification
of the Cauchy-Schwartz argument in the beginning of the proof of Lemma~\ref{A.8}
yields
$$
|\widetilde L_\sigma^{mn _3} (g)(x)|^2
\le K \frac {\lambda_{2 \sigma}^{mn_3} }{ \lambda_\sigma^{2mn_3}} 
\widetilde L_0^{mn _3} |g^2|(x)\, ,
$$
for some  $K$ independent of $mn_3$ and $f$.
Next, assume $\sup |f'| < 2 C |t| \sup |f| $ and
assume $\|f \|_{1,t}= 1$. By the spectral properties of
$\widetilde L_0$ on the space of Lipschitz functions endowed
with the norm $\sup |g| + \Lip(g)$ (with
$\Lip(g)$ the smallest Lipschitz constant of $g$), there are
$R_\sigma<\infty$,  $\tau^{L}_\sigma <1$ (independent of $f$ and $t$),
with:
\begin{equation*}
\begin{split}
\sup& |\widetilde L_s^{2mn_3}(f)|^2\le 
\sup |\widetilde L_\sigma^{mn_3} (\widetilde L_s^{mn_3}(f))|^2
\le  K \frac {\lambda_ {2 \sigma}^{mn_3} }{ \lambda_\sigma^{2mn_3}} 
\sup \widetilde L_0^{mn _3} (|\widetilde L_s^{mn_3}(f)|^2) \cr
&\quad\le K \frac {\lambda_ {2 \sigma}^{mn_3} }{ \lambda_\sigma^{2mn_3}}
 \bigl ( \int |\widetilde L_s^{mn_3}(f)|^2\, d\nu_0
+ R_\sigma (\tau^L_\sigma)^{mn_3} [\sup+ \Lip](|\widetilde L_s^{2mn_3}(f)|^2)\bigr ) \cr
&\quad \le K \frac {\lambda_ {2 \sigma-1}^{mn_3} }{ \lambda_\sigma^{2mn_3}}
\cdot \biggl (\sup |f|^2 \beta^{m}+
 R_\sigma (\tau^L_\sigma)^{mn_3}  [\sup+\Lip](|\widetilde L_s^{mn_3}(f) |^2)\biggr ) 
\end{split}
\end{equation*}
using Lemma \ref{A.8} for $n=n_3$ and Cauchy-Schwartz). Lemma~\ref{A.7} gives
\begin{equation}
\label{tA.25}
\sup|\widetilde L_s^{mn_3}(f)|^2=\sup |f|^2 |v_m|^2\le \sup |f|^2 u_m^2\le \sup |f|^2\le 1 \, ,
\end{equation}
and
$\Lip(\widetilde L_s^{mn_3}(f)^2)\le 
2 \sup |f|^2 \cdot\sup |f|\sup |v'_m|\le 2 \sup |f|^3 2C |t|  \le 4 C |t|$, 
since $\Lip(|v_m|)\le \Lip(v_m)= \sup |v_m'|$.

In order to find $\max(\beta, \tau_\sigma^{n_3})
\cdot \frac {\lambda_ {2 \sigma-1}^{n_3} }{ \lambda_\sigma^{n_3}}
<\gamma_2 ^2<1$ so that (for all $m$)
$$
 K \frac {\lambda_ {2 \sigma-1}^{mn_3} }{ \lambda_\sigma^{2mn_3}}
\cdot \biggl ( \beta^{m}+  R_\sigma \tau_\sigma^{mn_3} (1+4C|t|)\biggr ) 
\le   \gamma_2^{2m} \, ,
$$
it is enough to require $m \ge \widehat A \log|t|$ for some  $\widehat A>0$
(and $\sigma$ close enough to $0$).

\noindent To control the derivative, invoke
Lemma~\ref{A.3}, exploiting the bounds just obtained:
\begin{equation*}
\begin{split}
\frac {\sup |(\widetilde L_s^{2mn_3}(f))'|}{ |t|}
&\le \frac { C(\Sigma,\bar K) |s|}{ |t|} \sup (\widetilde L_\sigma^{2mn_3}|f|) +
\frac {\rho^{mn_3} }{ |t|} \sup (\widetilde L_\sigma^{2mn_3}|(f)'|)\cr
&\le  \frac { C(\Sigma,\bar K) |s|}{ |t|}\gamma_2^m +
2 C \rho^{mn_3} \sqrt K \frac {\lambda_ {2 \sigma-1}^{mn_3} }{ \lambda_\sigma^{mn_3}}  
\le  \gamma_3^m  \, .
\end{split}
\end{equation*}
Take $n_4=2n_3$ and large enough $A \ge \widehat A$.
\end{proof}

\bibliographystyle{amsplain}

\end{document}